\def\Bbb#1{{\bf #1}}

\def\fnote#1{\footnote}
\def\blacksquare{\hbox{\vrule width 4pt height 4pt depth 0pt}}

\def\cwleftpar#1#2{\leftskip #1 \rightskip #2 plus 1fill}
\def\cwrightpar#1#2{\leftskip #1 plus 1fill \rightskip #2}
\def\cwcenterpar#1#2{\leftskip #1 plus 1fill \rightskip #2 plus 1fill}
\def\cwfullpar#1#2{\leftskip#1\rightskip#2}

\def\cwoutdent#1#2{\llap{\hbox to #1{#2 \hss}}\ignorespaces}
\def\cwparbegin#1#2#3#4#5{
	\ifcase #1 \cwleftpar{#2}{#3}
	\or \cwrightpar{#2}{#3}
	\or \cwcenterpar{#2}{#3}
	\else \cwfullpar{#2}{#3}\fi
	\ifcase #4 \baselineskip = 1.5\baselineskip
	\or \baselineskip = 2\baselineskip
	\or \baselineskip = 3\baselineskip
	\else \baselineskip = 1\baselineskip\fi
	\ifdim #5 > 0in \else \noindent \fi
	\noindent\ignorespaces}
\documentclass{article}
\begin{document}
\advance \vsize by -1\baselineskip
\def\makefootline{
{\vskip \baselineskip \noindent \folio                                  \par
}}

\vspace*{2ex}
\noindent {\Huge Linear Transports along Paths in\\[0.5ex] Vector
							Bundles}\\[1.8ex]
\noindent {\Large III. Curvature and Torsion}
\vspace*{2ex}

\noindent Bozhidar Zakhariev Iliev
\fnote{0}{\noindent $^{\hbox{}}$Permanent address:
Laboratory of Mathematical Modeling in Physics,
Institute for Nuclear Research and \mbox{Nuclear} Energy,
Bulgarian Academy of Sciences,
Boul.\ Tzarigradsko chauss\'ee~72, 1784 Sofia, Bulgaria\\
\indent E-mail address: bozho@inrne.bas.bg\\
\indent URL: http://theo.inrne.bas.bg/$\sim$bozho/}

\vspace*{2ex}

{\bf \noindent Published: Communication JINR, E5-93-261, Dubna, 1993}\\[1ex]
\hphantom{\bf Published: }
http://www.arXiv.org e-Print archive No.~math.DG/0502008\\[2ex]

\noindent
2000 MSC numbers: 53C99, 53B99, 57R35\\
2003 PACS numbers: 02.40.Ma, 02.40.Vh, 04.90.+e\\[2ex]

\noindent
{\small
The \LaTeXe\ source file of this paper was produced by converting a
ChiWriter 3.16 source file into
ChiWriter 4.0 file and then converting the latter file into a
\LaTeX\ 2.09 source file, which was manually edited for correcting numerous
errors and for improving the appearance of the text.  As a result of this
procedure, some errors in the text may exist.
}\\[2ex]

	\begin{abstract}
 Curvature and torsion of linear transports along paths in, respectively,
vector bundles and the tangent bundle to a differentiable manifold are
defined and certain their properties are derived.
	\end{abstract}\vspace{3ex}

 {\bf 1. INTRODUCTION}

\medskip
This work is devoted to the introduction of the concepts "curvature" and "torsion" with respect to the linear transports along paths considered in [1]. The former concept turns out to be more general in a sense that it can be introduced in arbitrary vector bundles, while the latter one can be defined only in the tangent bundle to a differentiable manifold.

The theory of linear transports along paths will not be repeated here, the reader being referred to [1] for further details. We shall only mention briefly hereafter the elements which are necessary for the present part of our investigation.

By $(E,\pi ,B)$ we denote a vector bundle with a base $B$, total bundle space $E$ and projection $\pi :E  \to  $B. A path in $B$ is a map $\gamma :J  \to  B, J$ being a real interval (of arbitrary type).

To any $C^{1}$ linear transport $(L$-transport) along paths $L$ in $(E,\pi
,B)$ is associated a derivation ${\cal D}$ along paths which maps $\gamma $
into a derivation ${\cal D}^{\gamma }$along $\gamma $ such that its action on
any $C^{1}$ section $\sigma $ over $\gamma (J)$ is $(s$ and $s+\epsilon $
belong to $J)$
\[
 {\bigl(}{\cal D}^{\gamma }\sigma {\bigr)}(\gamma (s)):={\cal
D}^{\gamma }_{s}\sigma
:=\lim_{\varepsilon\to0}
\bigl[\frac{1}{\varepsilon} {\bigl(}
L^{\gamma}_{s+\epsilon \to  s}\sigma (\gamma (s+\epsilon ))
-\sigma (\gamma (s)){\bigr)} \bigr], \qquad (1.1)
\]
 where $L^{\gamma}_{s\to t}, s,t\in J$ is the $L$-transport along $\gamma $
from $s$ to t.

If $\{e_{i}\} ($the Latin indices run from 1 to $\dim(\pi ^{-1}(x)), x\in
B)$ is a field of bases along $\gamma $, i.e. $\{e_{i}(s)\}$ is a basis in
$\pi ^{-1}(\gamma (s))$, and $\sigma =\sigma ^{i}e_{i}(a$ summation from 1 to
$\dim(\pi ^{-1}(x)), x\in B$ over repeated on different levels indices is
assumed), then the explicit form of (1.1) is
\[
{\bigl(}{\cal D}^{\gamma }\sigma {\bigr)}(\gamma (s))
= \Bigl[ \frac{d\sigma^i(\gamma(s))}{ds}
+ \Gamma ^{i}_{.j}(s;\gamma )\sigma ^{j}(\gamma (s))
\Bigr]  e_{i}(s)\qquad (1.2)
\]
in which $\Gamma ^{i}_{.j}(s;\gamma )$ are the coefficients of $L
($coinciding with that of ${\cal D})$. The general form of $\Gamma
^{i}_{.j}(s;\gamma )$ is given by
\[
\Gamma _{\gamma }(s)= [\Gamma ^{i}_{.j}(s;\gamma )]
  =F^{-1}(s;\gamma ) \frac{dF(s;\gamma)}{ds}
\quad
\Bigl( =
\frac{\partial H(s,t;\gamma)}{\partial t}\Big|_{t=s }
\Bigr) ,\qquad (1.3)
\]
 $F(s;\gamma )$ being nondegenerate
$C^{1}$matrix function of $s ($defining also the general form of the matrix
$H(s,t;\gamma )$ of $L; cf.  [1])$. If along $\gamma $ is made the change
$\{e_{i}(s)\} \to \{e_{i^\prime }(s)=A^{i}_{i^\prime }(s)e_{i}(s)\}$ of the
basis $\{e_{i}\}$, where
$[ A^{i}_{i^\prime }(s)] =: [ A^{i^\prime }_{i}(s)]^{-1}$
is a nondegenerate matrix function, then in the new basis
$\{e_{i^\prime }\}$ the coefficients of $L$ are
\[
\Gamma ^{i^\prime }_{..j^\prime }(s;\gamma )=A^{i^\prime
}_{i}(s)A^{j}_{j^\prime }(s)\Gamma ^{i}_{.j}(s;\gamma )
+A^{i^\prime }_{i}(s) \frac{d A_{j'}^{i}(s)}{ds}.  \qquad (1.4)
\]
If $B$ is a manifold and $\gamma :J  \to  B$ is a $C^{1}$ path with a tangent
vector field, then of special interest are the $L$-transports along
$C^{1}$ paths whose coefficients have the form
\[
 \Gamma ^{i}_{.j}(s;\gamma )
=\sum_{\alpha=1}^{\dim(B)} \Gamma ^{i}_{.j\alpha }(\gamma(s))
\dot\gamma^{\alpha }(s) \qquad (1.5)
\]
for some functions $\Gamma ^{i}_{.j\alpha }:B  \to  {\Bbb R}$. In particular,
if $(E,\pi ,B)$ is the tangent bundle to $B$, then $\Gamma ^{i}_{.j\alpha
}$are coefficients of a linear connection, the corresponding to which
parallel transport coincides with the considered $L$-transport along (smooth)
paths.

The usage of the operator ${\cal D}^{\gamma }$, associated to a given $L$-transport along paths through (1.1), gives a possibility for defining the torsion and the curvature of that transport to which is aimed this paper. The torsion is introduced in Sect. 2, where it is prove that it vanishes only for parallel transports generated by symmetric linear connections. The curvature is defined in Sect. 3, where also some its properties are investigate. Certain concluding remarks are made in Sect. 4.

\medskip
\medskip
 {\bf 2. TORSION}

\medskip
Let $M$ be a differentiable manifold and $(T(M),\pi ,M)$ be its tangent
bundle (some times denoted simply by $T(M)) [2]$.

Let $\eta :J\times  J^\prime   \to  M, J$ and $J^\prime $ being real
intervals, be a $C^{2}$map. With $\eta (\cdot ,t)$ and $\eta (s,\cdot ),
s,t\in J$ we denote the paths $\eta (\cdot ,t):s\to   \to  \eta (s,t)$ and
$\eta (s,\cdot ):t\to  \eta (s,t)$ and with $\eta ^\prime (\cdot ,t)    \eta
^{\prime\prime}(s,\cdot )$, respectively, the tangent to them vector fields.

{\bf Definition 2.1.} The torsion operator of an $L$-transport along paths
in $(T(M),\pi ,M)$ is the map
${\cal T}:\eta \to  {\cal T}^{\eta }:J\times J^\prime \to T(M)$
defined for all $(s,t)\in J \times J^\prime $ by
\[
{\cal T}^{\eta }(s,t):={\cal D}^{\eta (\cdot ,t)}_{s}\eta
^{\prime\prime}(\cdot ,t)-{\cal D}^{\eta (s,\cdot )}_{t}\eta ^\prime (s,\cdot
)\in T_{\eta (s,t)}(M),\qquad (2.1)
\]
 where ${\cal D}$ is the derivation along
paths corresponding in accordance with (1.1) to the given $L$-transport along
paths.

 Let in $(T(M),\pi ,M)$ be fixed some local coordinate basis. For the sake of
simplicity we assume $\eta (J,J^\prime )$ to lie in only one coordinate
neighborhood of the generating this basis coordinates. Then, using (1.2), we
get the components of ${\cal T}^{\eta }(s,t)$ to be
\[
 ( {\cal T}^{\eta }(s,t) )^{ i}
=
\Bigl( \frac{\partial}{\partial s} \frac{\partial}{\partial t}\eta ^{i  }
\Biggr)\Big|_{(s,t)}
 + \Gamma^{i}_{.j}(s;\eta (\cdot ,t))(\eta^{\prime\prime}(s,t))^{j} -
\]
\[
 -
\Bigl( \frac{\partial}{\partial t}\frac{\partial}{\partial s}\eta ^{i  }
\Biggr)\Big|_{(s,t)}
+ \Gamma ^{i}_{.j}(t;\eta (s,\cdot ))(\eta ^\prime (s,t))^{j},
\]
 where $\Gamma ^{i}_{.j}(\ldots  )$ are the coefficients of the considered
$L$-transport. As we supposed $\eta $ to be $a C^{2}$map, from here we find
\[
 {\cal T}^{\eta }(s,t)^{  i}_{  }=\Gamma ^{i}_{.j}(s;\eta (\cdot ,t))(\eta
^{\prime\prime}(s,t))^{j}-\Gamma ^{i}_{.j}(t;\eta (s,\cdot ))(\eta ^\prime
(s,t))^{j}.  \qquad (2.2)\]

{\bf Proposition 2.1.} If an $L$-transport in  $(T(M),\pi ,M)$ is a parallel
transport generated by a linear connection $\nabla $ with local coefficients
$\Gamma ^{i}_{.jk}$, then in every coordinate basis
\[
\bigl( {\cal T}^{\eta }(s,t) \bigr)^{  i}
=[\Gamma ^{i}_{.jk}(\eta (s,t))-\Gamma
^{i}_{.kj}(\eta (s,t))](\eta ^{\prime\prime}(s,t))^{j}(\eta ^\prime
(s,t))^{k}.  \qquad (2.3)
\]

 {\bf Proof.} From what is supposed it follows that the
coefficients of the $L$-transport have the representation (1.5), i.e.
\[
 \Gamma ^{i}_{.j}(s;\gamma )
= \Gamma ^{i}_{.jk}(\gamma (s)) \dot\gamma^{k}(s)\qquad (2.4)
\]
for every  $\gamma $.  Substituting this equality into (2.2), we get
(2.3).\blacksquare

Under the condition of the above proved proposition it naturally arises the
torsion tensor $T$ of $\nabla $. In fact, if we take two $C^{2}$vector fields
A and $B$ defined on one and the same neighborhood $U\subset M$ and for every
$x\in U$ define
$\eta _{x}$ by
$\eta _{x}(s_{0},t_{0})=x,$
$\eta'_x(s,t_{0})=A_{\eta_{x}(s,t_0)}$ and
$\eta^{\prime\prime}_x(s_{0},t)=B_{\eta_{x}(s_0,t)}$ for a fixed
$(s_{0},t_{0})\in J\times  J^\prime $, then, we can put
\[
(T(A,B))(x):={\cal T}^{\eta _{x}}(s_{0},t_{0}).\qquad (2.5a)
\]
 So, from (2.3), we obtain
\[
 [ T(A,B))(x) ]^{i}=[-\Gamma ^{i}_{.jk}(x) + \Gamma
^{i}_{.kj}(x)]A^{j}(x)B^{k}(x),\qquad (2.5b)
\]
 hence $T(A,B)=\nabla
_{A}B-\nabla _{B}A-[A,B], [A,B]$ being the commutator of A and $B [2]$, is
the usual torsion tensor $(cf. [2])$.

The next proposition shows that the torsion operator, as in a case of
parallel transports generated by linear connections, characterizes the
deviation in the general case of arbitrary $L$-transports in the tangent
bundle to a given manifold from the parallel transports generated by
symmetric affine connections.

{\bf Proposition 2.2.} The torsion operator of a given $L$-transport in
$(T(M),\pi ,M)$ is zero in a neighborhood of some point $x\in M$, i.e.
\[
 {\cal T}^{\eta }(s,t)=0\quad for\ every \ \eta :J\times J^\prime \to  M,
\ \eta (s_{0},t_{0})=x,\qquad  \qquad (2.6)
\]
 if and only if that transport is a
 parallel transport generated by symmetric linear connection in that
neighborhood.

{\bf Proof.} The sufficiency is almost evident: if the given $L$-transport is a parallel transport generated by a symmetric linear connection, then $T=0 ($by definition of such connection) and it is valid (2.4) as a consequence of which from (2.5a) follows (2.6).

 Conversely, let (2.6) be true. Then from  (2.2), we get
\[
\Gamma ^{i}_{.j}(s;\eta (\cdot ,t))
(\eta ^{\prime\prime}(s,t))^{j}
=\Gamma ^{i}_{.l}(t;\eta (s,\cdot ))(\eta ^\prime (s,t))
\]
from which follows
\[
 \Gamma ^{i}_{.j}(r;\gamma )=\Gamma ^{i}_{.jl}(r;\gamma ) \dot\gamma^{l}(r),
\]
 where $\Gamma ^{i}_{.jl}(r;\gamma )$ are some functions, $\gamma =\eta
(\cdot ,t),\eta (s,\cdot )$, respectively $r=s,t$ and    is the tangent to
$\gamma $ vector field.

Therefore, we see that $\Gamma ^{i}_{.jk}(s;\eta (\cdot ,t))=\Gamma
^{i}_{.kj}(t;\eta (s,\cdot ))$. So, as a consequence of the arbitrariness of
$\eta $, the common value of the right and left hand sides of this equality
can depend only on the point $\eta (s,t)$. Denoting this common value by
$\Gamma ^{i}_{.kj}(\eta (s,t))$ we obtain $\Gamma ^{i}_{.kl}(\eta
(s,t))=\Gamma ^{i}_{.lk}(\eta (s,t))$ and
\[
 \Gamma^{i}_{.k}(r;\gamma )
=\Gamma^{i}_{.kl}(\gamma (r)) \dot\gamma^{l}(r),
\quad r\in J,\  \gamma =\eta (\cdot ,t),\eta (s,\cdot ),
\quad resp.\ r=s,t.
\]

 As $\eta $ is arbitrary, from here we find
\[
\Gamma ^{i}_{.k}(r;\gamma )=\Gamma ^{i}_{.kl}(\gamma (r))\dot\gamma^{l}(r),
\quad \Gamma^{i}_{.kl}=\Gamma ^{i}_{.lk}\qquad (2.7)
\]
 for every path $\gamma :J  \to  M$
and $r\in $J. Due to (2.5b) and proposition 5.1 of [1] (see also the comments
after its proof) this means that $\Gamma ^{i}_{.kl}$are coefficients of a
symmetric linear connection the generated by which parallel transport
coincides with the initial $L$-transport (see also proposition 5.4 from [3]
and the note after its proof).\blacksquare

\medskip
\medskip
 {\bf 3. CURVATURE}

\medskip
The curvature operator, which is referred to an arbitrary vector bundle $(E,\pi ,B)$ is defined analogously to the torsion operator.

 Let $\eta :J\times  J^\prime   \to  $B.

{\bf Definition 3.1.} The curvature operator of an $L$-transport along paths
in the vector bundle $(E,\pi ,b)$ is a map ${\cal R}:\eta \to  {\cal R}^{\eta
}, {\cal R}^{\eta }:(s,t)  \to  {\cal R}^{\eta }(s,t)$ where ${\cal R}^{\eta
}(s,t)$ maps $C^{2}$sections of $(E,\pi ,b)$ into sections of $(E,\pi ,b)$
and it is defined for every $(s,t)\in J\times  J^\prime $ by
\[
{\cal R}^{\eta }(s,t):={\cal D}^{\eta (\cdot ,t)}\circ {\cal D}^{\eta
(s,\cdot )}- {\cal D}^{\eta (s,\cdot )}\circ {\cal D}^{\eta (\cdot
,t)}.\qquad (3.1)
\]
 If A is a $C^{2}$ section of $(E,\pi ,B)$ and $\{e_{i}\}$
is a field of bases in $\eta (J,J^\prime )$, i.e. $\{e_{i}(s,t)\}$ is a basis
in $\pi ^{-1}(\eta (s,t))$, then, by means of (1.2), we get
\[
 [({\cal R}^{\eta }(s,t))(A)](\eta (s,t))={\cal R}^{\eta }(s,t)^{  i}_{
.j}A^{j}(\eta (s,t))e_{i}(s,t),\qquad (3.2)
\]
 where
\[
 {\cal R}^{\eta }(s,t)^{  i}_{  .j}
= \frac{\partial }{\partial s}
	\bigl( \Gamma ^{i}_{.j}(t;\eta(s,\cdot ))\bigr)
 -   \frac{\partial }{\partial t}
      \bigl( \Gamma ^{i}_{.j}(s;\eta (\cdot ,t))\bigr) +
\]
\[
+ \Gamma ^{i}_{.k}(s;\eta (\cdot ,t))\Gamma ^{k}_{.j}(t;\eta (s,\cdot
))-\Gamma ^{i}_{.k}(t;\eta (s,\cdot ))\Gamma ^{k}_{.j}(s;\eta (\cdot
,t))\qquad (3.3)
\]
 are the components of ${\cal R}$ with respect to
$\{e_{i}\}$ at the point $(\eta ,s,t)$.

{\bf Proposition 3.1.} If $B$ is a manifold, $\eta:J\times J^\prime \to B$ is
of class $C^{2}$and the coefficients of some $L$-transport have the form
(1.5), then the components of its curvature operator in any basis are
\[
  {\cal R}^{\eta }(s,t)^{  i}_{  .j}={ } R(\eta (s,t))^{  i}_{  .j\alpha
\beta }(\eta ^\prime (s,t))^{\alpha }(\eta ^{\prime\prime}(s,t))^{\beta
},\qquad (3.4)
\]
 where
\[
 R(\eta (s,t))^{  i}_{  .j\alpha \beta }
= \Bigl[- \frac{\partial }{\partial x^\beta}
	(\Gamma ^{i}_{.j\alpha })
	+\frac{\partial }{\partial x^\alpha}
	(\Gamma ^{i}_{.j\beta })
-\Gamma ^{k}_{.j\alpha }\Gamma ^{i}_{.k\beta }
+\Gamma ^{k}_{.j\beta }\Gamma ^{i}_{.k\alpha   }
\Bigr]
\Big|_{\eta (s,t)}. (3.5)
\]

{\bf Proof.} This proof reduces to a simple substitution of (1.5) into
(3.3).\blacksquare

 If an $L$-transport is a parallel transport in $(T(M),\pi ,M)$ generated by
a linear connection $\nabla $ with coefficients $\Gamma ^{i}_{.jk}$, then, as
a consequence of the last proposition, naturally arises the curvature tensor
$R$ of $\nabla  [2]$. Actually, let us define $A, B$ and $\eta _{x}$as it was
done before (2.5a) and put
\[
  (R(A,B)C)(x):=[({\cal R}^{\eta _{x}}(s_{0},t_{o}))C](x)\qquad (3.6)
\]
 for
arbitrary vector field C. From here, (3.4) and (3.5), we find $R(A,B)=\nabla
_{A}\nabla _{B}-\nabla _{B}\nabla _{A}-\nabla _{[A,B]}$and that the
components of $R$ at $x\in M$ are exactly $(R(x))^{i}_{\hbox{.jkl}}$as they
are given by $(3.5) (cf$. e.g. [2]).

 {\bf Proposition 3.2.} Let in some basis $\{e_{i^\prime }\}$ the
coefficients of an $L$-transport in a vector bundle $(E,\pi ,B)$ be zeros
along every path, i.e.
\[
 \Gamma ^{i^\prime }_{..j^\prime }(s;\gamma )=0
\qquad for\ every\ \gamma :J  \to B. \qquad (3.7)
\]
Then:

1) If $B$ is a manifold and $\gamma $ is $a C^{1}$path, then in every basis
$\{e_{i}\}$ is valid (1.5), i.e.
\[
  \Gamma ^{i}_{.j}(s;\gamma )=\Gamma ^{i}_{.j\alpha }(\gamma (s))
\dot\gamma^{\alpha}(s), \gamma :J  \to  B,\qquad (3.8)
\]
 for some functions $\Gamma
^{i}_{.j\alpha }:B  \to  {\Bbb R}$. Besides, in the basis $\{e_{i^\prime }\}$
is fulfilled the equality
\[
 \Gamma ^{i^\prime }_{..j^\prime \alpha }(x)=0, x\in B;\qquad (3.9)
\]

 2) The curvature of the considered $L$-transport is zero, i.e.
\[
  {\cal R}^{\eta }(s,t)=0\qquad
for\ every\ \eta :J\times  J^\prime \to  B,\qquad (3.10)
\]
 which, if $B$ is a manifold, is equivalent to
\[
R(x)^{  i}_{  .j\alpha \beta }=0,\quad  x\in B. \qquad (3.10^\prime )
\]

{\bf Proof.} Let us choose a basis $\{e_{i}\}$ in $E$ such that $e_{i}(x)=
=A^{i^\prime }_{i}(x)e_{i^\prime }(x), x\in B, \det[A^{i^\prime }_{i}] \neq
0,\infty $ and $  [A^{i}_{i^\prime }] := [ A^{i^\prime }_{i}]  ^{-1}$.

 	From (3.7) and (1.4) follows
$0=\Gamma ^{i}_{.j}(s;\gamma )A^{j}_{j^\prime }(\gamma (s))
+ \frac{d}{ds}A^{j}_{j^\prime }(\gamma (s))$ or
\[
\Gamma^{i}_{.j}(s;\gamma )
=- \Bigl( \frac{d}{ds} A^{i}_{j^\prime }(\gamma (s))\Bigr))
 A^{j^\prime }_{j}(\gamma (s))
=\Bigl(\frac{\partial A_{j'}^{\ i}}{\partial s} A^{j^\prime }_{j}(x)
 \Bigr)\Big|_{ x=\gamma(s)} \dot\gamma^{\alpha }(s) \qquad (3.11)
\]
 where
$\{x^{\alpha }\}$ are local coordinates in a neighborhood of $\gamma (s)$.
Putting
\[
 \Gamma ^{i}_{.j\alpha }(x)
=- \frac{\partial A_{j'}^{i}}{\partial  x^\alpha} A^{j^\prime }_{j}(x),
\qquad (3.12)
\]
 we see from here that in any basis $\{e_{i}\}$ is true (3.8).

 In particular, if we let $e_{i}=e_{i^\prime }$, i.e. $A^{i^\prime
}_{i}(x)=\delta ^{i^\prime }_{i}$, we find
$\Gamma ^{i^\prime }_{..j^\prime \alpha }(\gamma (s))\dot\gamma^{\alpha
}(s)=\Gamma ^{i^\prime }_{.,j^\prime \alpha }(s;\gamma )=0$ for every path
$\gamma $, from which immediately follows (3.9).

At the end, the equality (3.10) is a corollary from (3.7) and (3.3), written in $\{e_{i^\prime }\}$, and $(3.10^\prime )$ is valid because of proposition 3.1, written in the basis $\{e_{i^\prime }\}$, in which is true the already proved equality (3.9).\blacksquare

{\bf Proposition 3.3.} If the curvature operator of an $L$-transport in a
vector bundle $(E,\pi ,B)$ is zero, i.e.
\[
{\cal R}^{\eta }(s,t)=0\quad
for\ every\ \eta :J\times J^\prime  \to  B,\ s,t\in J,\qquad (3.13)
\]
 then in $E$ there exist  fields of bases $\{e_{i^\prime }\}$
in which the coefficients of the $L$-transport are zeros along every path,
i.e. (3.7) is valid. Besides, all of these field of bases are obtained from
one another by linear transformations with constant coefficients.

{\bf Proof.} Let us choose a field of bases $\{e_{i}\}$ in E. We have to
prove the existence of a transformation $\{e_{i}(x)\}  \to  \{e_{i^\prime
}(x)=A^{i}_{i^\prime }(x)     e_{i}(x), \det[ A^{i}_{i^\prime }(x)] \neq
0,\infty \}$, such that in the basis $\{e_{i^\prime }\}$ to be fulfilled
(3.7).

If we let $[A^{i^\prime }_{i}(x)]:= [ A^{i}_{i^\prime }(x)] ^{-1}$, then,
due to (1.4), it follows that (3.7) is equivalent to
\[
\frac{d}{ds} A^{i^\prime }_{i}((\gamma (s))=A^{i^\prime }_{l}(\gamma
(s))\Gamma ^{l}_{.i}(s;\gamma )
\qquad for\ every\ \gamma :J  \to  B.      \qquad (3.14)
\]
Now we shall prove the existence of $A^{i^\prime }_{i}(x)$ satisfying this
equality which will mean also the existence of the pointed transformation.

Let $V:=J\times\cdot\times J (n=\dim(B)$ times) and $\eta :V  \to  $B.
For every $y=(y_{1},\ldots  ,y_{n})\in V$ we define maps
$\eta ^{y}_{i}:J \to  B$ and $\eta ^{y}_{ij}:J\times J^\prime \to  B$,
such that $\eta^{y}_{i}(s):=\eta (y)|_{y_{i}=s}$ and
$\eta^{y}_{ij}(s,t):=\eta (y)|_{y_{i}=s,y_j=t}$ (for $i=j$ we must have $s=t$
here).  Evidently $\eta ^{y}_{i}$ does not depend on $y_{i}$and $\eta
^{y}_{ij}$is independent of $y_{i}$ and $y_{j}$.

 Now we shall prove the existence of $A^{i^\prime }_{i}(x)$ satisfying
\[
\frac{\partial}{\partial y^j} A^{i^\prime }_{i}(\eta (y))
= A^{i^\prime }_{l}(\eta (y))\Gamma ^{l}_{.i}(y_{j};\eta ^{y}_{j}),
\qquad (3.15)
\]
 which is obtained from (3.14) for $s=y_{j}$and $\gamma =\eta _{j}$, as a
consequence of which $\eta _{j}(y_{j})=\eta (y)$.

	 In fact, a simple calculation shows that
\[
\Bigl(
\frac{\partial^2}{\partial y^j \partial y^k} -
\frac{\partial^2}{\partial y^k \partial y^j}
\Bigr)
A^{i^\prime }_{i}(\eta (y))
= \frac{\partial}{\partial y^j}
\Bigl(A^{i^\prime }_{l}(\eta (y))\Gamma ^{l}_{.i}(y_{k};\eta ^{y}_{k})
\Bigr)-
\]
\[
- \frac{\partial}{\partial y^k}
\Bigl(
A^{i^\prime }_{l}(\eta (y))\Gamma ^{l}_{.i}(y_{j};\eta ^{y}_{j})
\Bigr)
 = A^{i^\prime }_{m}(\eta (y))
\bigl( {\cal R}^{\eta ^{y}_{jk}}(y_{j},y_{k}) \bigr)^{  m}_{ .i},
\]
 from which, due to (3.13), follows the expression between the first and the
second equality signs to be equal to zero. But these are exactly the
integrability conditions for (3.15). Hence, the system (3.15) has a solution
with respect to $A^{i^\prime }_{i}. ($All such solutions are obtained from
one of them by a left multiplication with a constant matrix, but this is
insignificant here.)

It is important to be noted that by its construction this solution depends only on the point $\eta (y)$ but not on the map $\eta $ itself. Consequently, if we take a path $\gamma :J  \to  B$ and put $y^{j}=s, \eta _{j}=\gamma $ and $\eta (y)=\gamma (s)$ in (3.15) (the remaining arbitrariness in $\eta $ is insignificant), we get (3.14). So, in the basis $\{e_{i^\prime }\}$ defined by
$e_{i}(x)=A^{i^\prime }_{i}(x)e_{i^\prime }(x)$ is valid (3.7), i.e. in $\{e_{i^\prime }\}$ the coefficients of the considered $L$-transport are zeros.

At the end, if in two arbitrary bases $\{e_{i}\}$ and $\{e_{i^\prime
}=A^{i}_{i^\prime }e_{i}\}$ holds $\Gamma ^{i}_{.j}(s;\gamma )=\Gamma
^{i^\prime }_{..j^\prime }(s;\gamma )=0$ for every path $\gamma $, then from
(1.4) follows $dA^{i}_{i^\prime }(\gamma (s))/ds=0$, so we get $  A^{i^\prime
}_{i}(x)  =$const.\blacksquare

{\bf Theorem 3.1.} If in a vector bundle is given an $L$-transport along paths, then a necessary and sufficient condition for the existence of a field of bases in the bundle, in which the coefficients of the $L$-transport are zeros along any path, or, equivalently, in which the transport's matrix is unit along any path, is the curvature operator of this $L$-transport to be zero.

{\bf Proof.} This result is a direct consequence of propositions 3.3 and 3.2, as well as proposition 5.2 from [1].\blacksquare

{\bf Proposition 3.4.} If an $L$-transport in the tangent bundle to a manifold has a zero curvature, then a necessary and sufficient condition for the existence of a holonomic basis in which the coefficients of the $L$-transport are equal to zero is it to be torsion free, or, equivalently, when this $L$-transport is a parallel transport generated by a symmetric linear connection.

{\bf Proof.} First of all, the second part of the proposition is a corollary from proposition 2.2.

If the conditions of this proposition are fulfilled, then by proposition 3.2 in the pointed holonomic basis are valid $(3.7)-(3.9)$, the last of which, together with (2.3), gives ${\cal T}^{\eta }(s,t)=0$.

On the opposite, if ${\cal T}^{\eta }(s,t)=0$, then from the conditions of
the proposition, $(2.5a), (3.6)$, theorem 3.1 and proposition 3.2 follows
that the considered $L$-transport is a parallel transport generated by a
linear connection with vanishing curvature and torsion. Hence the
coefficients $\Gamma ^{i}_{.jk}$of this connection and $\Gamma ^{i}_{.j}$of
the $L$-transport are connected through the relation $\Gamma
^{i}_{.j}(s;\gamma )=\Gamma ^{i}_{.jk}(\gamma (s))\dot\gamma^{k}(s)$, but, as
is known (see e.g. [4]), for the first quantities there exists a holonomic
basis in which they vanish. In fact, if in $\{\partial /\partial x^{i}\}$ the
connection's coefficients are $\Gamma ^{i}_{.jk}$, then there exists $
A^{i^\prime }_{j}  =  A^{i}_{j^\prime }  $, such that $\Gamma
^{i}_{.jk}=A^{i}_{j^\prime }\partial A^{j^\prime }_{j}/\partial x^{k}$, as
the integrability conditions for these equations with respect to $A^{i^\prime
}_{j}$are exactly $R^{i}_{\hbox{.jkl}}=0$. The basis $e_{j^\prime
}=A^{i}_{j^\prime }\partial /\partial x^{i}$is holonomic because the
conditions for this are $T^{i}_{.jk}= =0$. And at the end, as a consequence
of (1.5), we have $\Gamma ^{i^\prime }_{..j^\prime k^\prime }=0.\blacksquare
$

\vspace{4ex}
 {\bf 4. CONCLUSION}

\medskip
In this work we have introduced the curvature and torsion of a linear transport along paths in, respectively, a vector bundle and the tangent bundle to a manifold. This was done by means of the assigned to the transport derivation along paths. Evidently, the same scheme works also with respect to any derivation along paths (we did not used explicitly the transports anywhere). This fact is in a natural agreement with the proved in [1] equivalence between $L$-transports along paths and derivations along paths. On the other hand, as the derivation along paths is a local concept such must be, and in fact are, the curvature and torsion of a linear transport along paths. From this point of view arises the problem for the global analogs of these concepts, which will be considered elsewhere.

At the end, we want to note that the investigated in [5] "flat" linear transports in tensor bundles are "flat" in the sense that their curvature operators, as they are defined in the present paper, are zero. This follows from the simple fact that the derivations along paths corresponding to these transports are covariant differentiations along paths generated by linear connections (see
[5], Sect. 3 and Sect. 3 of the present work).

\medskip
\medskip
 {\bf ACKNOWLEDGEMENTS}

\medskip
The author expresses his gratitude to Prof. Vl. Aleksandrov (Institute of Mathematics of Bulgarian  Academy of Sciences) for constant interest in this work and stimulating discussions.

This research was partially supported by the Fund for Scientific Research of Bulgaria under contract Grant No. $F 103$.

\medskip
\medskip
 {\bf REFERENCES}

\medskip
1.  Iliev B.Z., Linear transports along paths in vector bundles. I. General theory, Communication JINR, $E5-93-239$, Dubna, 1993.\par
2.  Kobayashi S., K. Nomizu, Foundations of Differential Geometry, Vol. 1, Interscience Publishers, New York-London, 1963.\par
3.  Iliev B.Z., Linear transports along paths in vector bundles. II. Some applications, Communication JINR, $E5-93-260$, Dubna, 1993.\par
4.  Schouten J. A., Ricci-Calculus: An Introduction to Tensor Analysis and its Geometrical Applications, Springer Verlag, Berlin-G\"ottingen-Heidelberg, $1954, 2-nd ed$.\par
5.  Iliev B.Z., Flat linear connections in terms of flat linear transports in tensor bundles, Communication JINR, $E5-92-544$, Dubna, 1992.

\newpage
 \noindent
\vspace*{5ex}
 Iliev B. Z. \\[3ex]

\noindent
 Linear Transports along Paths in Vector Bundles \\
 III. Curvature and Torsion \\[5ex]

 Curvature and torsion of linear transports along paths in, respectively,
vector bundles and the tangent bundle to a differentiable manifold are
defined and certain their properties are derived.\\[5ex]

 The investigation has been performed at the Laboratory of Theoretical
Physics, JINR.

\end{document}